\documentclass[11pt]{amsart}

\usepackage[utf8]{inputenc}  

\usepackage{amsmath}
\usepackage{amsfonts}
\usepackage{amssymb}
\usepackage{amsthm}
\usepackage{mathtools} 
\usepackage{latexsym}

\usepackage{enumerate}
\usepackage{cancel}
\usepackage{cases}
\usepackage{empheq}
\usepackage{multicol}

\usepackage{subfigure}
\usepackage{float}

\usepackage{wrapfig}

\usepackage{color}

\usepackage[backgroundcolor=gray!30,linecolor=black]{todonotes}

\usepackage{mathrsfs}
\usepackage{fontenc} 
\usepackage{inputenc} 

\usepackage{verbatim}

\theoremstyle{plain}
\newtheorem{theorem}{Theorem}[section]

\newtheorem{lemma}[theorem]{Lemma}

\newtheorem{remark}[theorem]{Remark}

\theoremstyle{definition}

\theoremstyle{remark}

\numberwithin{equation}{section} 
\numberwithin{figure}{section}   

\usepackage[square,comma,numbers,sort&compress]{natbib}
\usepackage{url}

\newcommand{\vect}[1]{\mathbf{#1}}

\newcommand{\bu}{\vect{u}}
\newcommand{\bv}{\vect{v}}
\newcommand{\bw}{\vect{w}}
\newcommand{\bA}{{A}}
\newcommand{\bB}{{B}}
\newcommand{\bP}{{P}}

\newcommand{\bx}{\vect{x}}

\newcommand{\maps}{\rightarrow}

\newcommand{\ip}[2]{\left<#1,#2\right>}

\newcounter{my_counter}
\title[Data assim. via Voigt with observable data]{
Approximate continuous data assimilation 
of the 2D Navier-Stokes equations
via the Voigt-regularization
with observable data}
\date{\today}

\author{Adam Larios}
\address[Adam Larios]{Department of Mathematics\\
                University of Nebraska--Lincoln\\
        Lincoln, NE 68588}
\email[Adam Larios]{alarios@unl.edu}
\author{Yuan Pei}
\address[Yuan Pei]{Department of Mathematics\\
                Western Washington University\\
        Bellingham, WA 98225}
\email[Yuan Pei]{yuan.pei@wwu.edu}

\thanks{MSC 2010 Classification: 
35A01, 
35B65, 
35K40, 
35K61,
35Q35, 
35Q30 
35Q93,
76D03. 
}

\thanks{\em Keywords:
Feedback control, 
Data assimilation, 
Navier-Stokes equations.
}

\begin{document}
\begin{abstract}
We propose a data assimilation algorithm for the 2D Navier-Stokes equations, based on the Azouani, Olson, and Titi (AOT) algorithm, but applied to the 2D Navier-Stokes-Voigt equations.  Adapting the AOT algorithm to regularized versions of Navier-Stokes has been done before, but the innovation of this work is to drive the assimilation equation with observational data, rather than data from a regularized system.  We first prove that this new system is globally well-posed.  Moreover, we prove that for any admissible initial data, the $L^2$ and $H^1$ norms of error are bounded by a constant times a power of the Voigt-regularization parameter $\alpha>0$, plus a term which decays exponentially fast in time.  In particular, the large-time error goes to zero algebraically as $\alpha$ goes to zero.  Assuming more smoothness on the initial data and forcing, we also prove similar results for the $H^2$ norm.
\end{abstract}

\maketitle
\thispagestyle{empty}

\noindent
\section{Introduction}\label{secInt}
\noindent
In real weather simulations, it is common practice to use some type of ``filtered'' or ``regularized'' model rather than the original equations to stabilize under-resolved simulations.  For instance, one might use a large-eddy simulation (LES) model or an $\alpha$-model, rather than the pure Navier-Stokes equations (see, e.g., \cite{Berselli_Iliescu_Layton_2006_book,Lesieur_Metais_Comte_2005_LES_book,Sagaut_2006_LES_book} and the references therein).  Somewhat separate from this approach is a class of techniques known as \textit{data assimilation}, which aim to improve the accuracy of simulations by incorporating observational data into the dynamical system being solved.  Combining these approaches is a natural idea, but a problem arises immediately because data from observations obviously does not originate from the regularized model, but from the underlying physical model\footnote{at least, in so far as the physical model can be considered to be correct.  A mismatch of physical model parameters is studied in \cite{CarlsonHudsonLarios_2018}.}.  This mismatch, and the error that results from it, is the subject of the present work.  In particular, we show that, in the setting we consider, the error decays exponentially fast in time, but only down to a level which is determined by the regularizing parameter $\alpha$.  However, we also show that this lower level (which results from the mismatch between models) goes to zero as a power of $\alpha$ when $\alpha\rightarrow0$.

We focus on the method of continuous data assimilation proposed by Azouani, Olson, and Titi in  \cite{Azouani_Olson_Titi_2014,Azouani_Titi_2014} (see also \cite{Cao_Kevrekidis_Titi_2001,Hayden_Olson_Titi_2011,Olson_Titi_2003} for early ideas in this direction), which we refer to as the AOT algorithm. For simplicity, we consider the 2D Navier-Stokes (NS) equations to be the underlying physical model, and take the ``filtered'' model to be the 2D Navier-Stokes-Voigt (NS-Voigt) equations.  While the 2D NS equations are already globally well-posed in the sense of possessing unique smooth solutions for all time, numerical simulations of these equations could still be under-resolved, and therefore, it is reasonable to consider filtered or smoothed versions of the equations.  We therefore consider an algorithm where data arises from the NS equations, but it is incorporated into the NS-Voigt equations via an AOT-type feedback-control term.  We show analytically that, so long as the parameters satisfy certain bounds, the error in the solutions decays exponentially fast in time until it reaches a certain level, and that moreover, this error level goes to zero as the Voigt-regularization parameter $\alpha$ goes to zero.

We note that the phenomenon of exponential decay of the error down to a level which can be controlled has recently been observed in a different context.  Namely, in \cite{CarlsonHudsonLarios_2018}, the effect of using mismatched viscosity terms is examined.  One viscosity (equivalently, Reynolds number) is used for the ``physical'' equation, and a different viscosity (or Reynolds number) is used for the ``assimilation'' equation.  Such a situation might arise, for example, when the true viscosity of the system is unknown, or when simulations use a smaller Reynolds number than the physical Reynolds number to stabilize the simulation.  It is proven in \cite{CarlsonHudsonLarios_2018} that the error decays exponentially down to a level that is controlled by the difference in the two viscosities, and that this level decreases to zeros as one viscosity approaches the other.  This result parallels the results of this paper, albeit in a very different context.

The NS-Voigt equations, and their inviscid counterpart, called the Euler-Voigt equations have been studied analytically and extended in a wide variety of contexts (see, e.g.,
 \cite{Berselli_Bisconti_2012,Berselli_Kim_Rebholz_2016,Berselli_Kim_Spirito_2016_DCDSB,Berselli_Spirito_2017,Bohm_1992,Borges_Ramos_2013,Cao_Lunasin_Titi_2006,Catania_2009,Catania_Secchi_2009,Cuff_Dunca_Manica_Rebholz_2015,Ebrahimi_Holst_Lunasin_2012,Gal_Medjo_2013_MMAS,Garcia_Luengo_Julia_Read_2012,Kalantarov_Levant_Titi_2009,Kalantarov_Titi_2009,Khouider_Titi_2008,Larios_Lunasin_Titi_2013,Larios_Lunasin_Titi_2015,Larios_Titi_2009,Larios_Titi_2010_MHD,Levant_Ramos_Titi_2009,Li_Qin_2013,Niche_2016_JDE,Olson_Titi_2007,Oskolkov_1973, Oskolkov_1982,Ramos_Titi_2010}, and the references therein).  Computational studies were carried out in \cite{Kuberry_Larios_Rebholz_Wilson_2012,DiMolfetta_Krstlulovic_Brachet_2015,Layton_Rebholz_2013_Voigt,Larios_Petersen_Titi_Wingate_2015}.  The NS-Voigt equations were first proposed by Oskolkov in \cite{Oskolkov_1973,Oskolkov_1982} as a model for Kelvin–Voigt fluids, but were later viewed as a regularization for the NS equations in \cite{Cao_Lunasin_Titi_2006}, where also the Euler-Voigt equations were first introduced and studied.  The Voigt-regularization is related to the wider class of  $\alpha$-models, including the NS-$\alpha$ (NS-$\alpha$) model and the Leray-$\alpha$ model, which were first proposed and studied in \cite{Cao_Lunasin_Titi_2006,Foias_Holm_Titi_2002,Ilyin_Lunasin_Titi_2006, Chen_Foias_Holm_Olson_Titi_Wynne_1998_PF,Chen_Foias_Holm_Olson_Titi_Wynne_1999,Holm_Titi_2005, Chen_Foias_Holm_Olson_Titi_Wynne_1998_PRL,Cheskidov_Holm_Olson_Titi_2005}.  
The Voigt model enjoys two features that the other $\alpha$-models are not known to have in the 3D case.  First, it is known to be globally well-posed in the inviscid case.  Second, in the viscous case, it is well-posed in the physical case of ``no-slip'' homogeneous Dirichlet boundary conditions, with no need to impose artificial boundary conditions.  Although we only work in 2D, we focus on the NS-Voigt model due in part to these attractive features, and also due to its simplicity.  The study of the AOT-algorithm applied to other models but still driven by observable data, and the errors resulting from the mismatch, will be the subject of a forthcoming work.

In \cite{Albanez_Nussenzveig_Lopes_Titi_2016}, it was shown that the AOT data assimilation scheme for the 3D NS-$\alpha$ model (with any admissible initial data) has solutions which converge to solutions of the 3D NS-$\alpha$ (with given admissible initial data).  A similar result was later proved in the context of the 3D-Leray-$\alpha$ model in \cite{Farhat_Lunasin_Titi_2018_Leray_AOT}.   We note that in both of these cases, the reconstruction of the NS-$\alpha$ solution is done via data that do not come from the NS equations themselves, but from the NS-$\alpha$ model.  In this sense, the observations are not physical, since the they are generated from NS-$\alpha$, which itself is not physical, but is used as a regularized version of Navier-Stokes; however, \cite{Albanez_Nussenzveig_Lopes_Titi_2016} and \cite{Farhat_Lunasin_Titi_2018_Leray_AOT} were important steps in demonstrating that AOT-type data assimilation is not limited to 2D equations, but can be extended to 3D in certain contexts.  That is to say, the barrier to extending theoretical results about the AOT method for 2D NS to 3D NS is not directly due to the dimensionality, but rather it is likely due to the nature of the equations in three dimensions, which is the same barrier standing in the way of resolving the major open problem of global existence for strong solutions of the 3D NS equations.  

It is worth noting that classical data assimilation is largely focused on statistical approaches utilizing the Kalman filter \cite{Kalman_1960_JBE} and 3D/4D-Var methods, and their many variations (see, e.g., \cite{Daley_1993_atmospheric_book,Kalnay_2003_DA_book,Law_Stuart_Zygalakis_2015_book,Lewis_Lakshmivarahan_2008}, and the references therein).  
The AOT algorithm, which is also referred to in the literature as \textit{continuous data assimilation}, differs markedly from the Kalman filter approach.  Rather than starting from the numerical level and employing statistical tools, AOT data assimilation manifests at the PDE level by way of a feedback-control term which penalizes deviations from interpolations of observable data.  The use of interpolation is a key difference between the AOT method and the so-called nudging or Newtonian relaxation methods introduced in \cite{Anthes_1974_JAS,Hoke_Anthes_1976_MWR}.  See, e.g., \cite{Lakshmivarahan_Lewis_2013} for an overview of nudging methods.  We note that a method that  resembles the AOT method in some respects was introduced in \cite{Blomker_Law_Stuart_Zygalakis_2013_NL} in the context of stochastic differential equations. 
Much recent literature has built upon the AOT algorithm and its ideas; see, e.g., \cite{Albanez_Nussenzveig_Lopes_Titi_2016,
Altaf_Titi_Knio_Zhao_Mc_Cabe_Hoteit_2015,
Bessaih_Olson_Titi_2015,
Biswas_Foias_Monaini_Titi_2018downscaling,
Biswas_Martinez_2017,
CarlsonHudsonLarios_2018,
Celik_Olson_Titi_2018,
Farhat_Jolly_Titi_2015,
Farhat_Lunasin_Titi_2016abridged,
Farhat_Lunasin_Titi_2016benard,
Farhat_Lunasin_Titi_2016_Charney,
Farhat_Lunasin_Titi_2017_Horizontal,
Foias_Mondaini_Titi_2016,
Foyash_Dzholli_Kravchenko_Titi_2014,
GarciaArchilla_Novo_Titi_2018,
Gesho_Olson_Titi_2015,
GlattHoltz_Kukavica_Vicol_2014,
Ibdah_Mondaini_Titi_2018uniform,
Jolly_Martinez_Olson_Titi_2018_blurred_SQG,
Jolly_Martinez_Titi_2017,
Jolly_Sadigov_Titi_2015,
Larios_Lunasin_Titi_2015,
Larios_Pei_2017_KSE_DA_NL,
Lunasin_Titi_2015,
Markowich_Titi_Trabelsi_2016,
Mondaini_Titi_2018_SIAM_NA,
Pei_2018,
Rebholz_Zerfas_2018_alg_nudge}.  
Computational trials of the AOT algorithm and its variants were carried out in the case of the NS equations \cite{Gesho_Olson_Titi_2015,Larios_Rebholz_Zerfas_2018,Leoni_Mazzino_Biferale_2018}, the 2D B\'enard convection equations \cite{Altaf_Titi_Knio_Zhao_Mc_Cabe_Hoteit_2015}, and the 1D Kuramoto-Sivashinsky equations \cite{Lunasin_Titi_2015,Larios_Pei_2017_KSE_DA_NL}. 



We now describe the main algorithm proposed in this work.  
First, consider the Navier-Stokes equations over the 2D periodic domain $\Omega = \mathbb{T}^2 = \mathbb{R}^2/\mathbb{Z}^2$,
\begin{subequations}\label{NSE}
 \begin{empheq}[left=\empheqlbrace]{align}
   \label{NSE1}
     \partial_t\bu 
     +
     (\bu\cdot\nabla)\bu
     +
     \nabla p
     &=
     \nu\Delta\bu 
     + 
     \mathbf{f},
     \\
   \label{NSE_div}
     \nabla\cdot \bu
     &= 0,
     \\
   \label{NSE0}
     \bu(\bx,0) &= \bu_0(\bx).
 \end{empheq}
\end{subequations}1
For the purposes of the present work, we take this system to be the correct underlying physical system.  In practice, the initial data $\bu_0$ might not be known at every point in the domain; indeed, it might not be known at all.  However, we assume that we have some information about $\bu(t)$, such as measurements\footnote{Here, for simplicity, we assume the measurements to be taken at every time $t$, and to have perfect accuracy, but see \cite{Celik_Olson_Titi_2018} for an adaption of the AOT algorithm to the case of discrete-in-time measurements, and see \cite{Bessaih_Olson_Titi_2015} for adaption to the case of stochastically noisy data.} of its values at certain points in the domain, from which we may construct an approximation of $\bu(t)$, namely $I_h(\bu(t))$, where $I_h$ is a linear operator based on $\mathcal{O}(h^{-1})$ data points, satisfying \eqref{Asml} below.  Thus, we will consider a new system with an AOT-type feedback-control term $I_h(\bu) - I_h(\bv)$.  Moreover, we assume that we have some reason to consider a regularized version of the original physical system; for example, it may be that our simulation is under-resolved, leading to numerical errors.  For the reasons described above, we use the Voigt regularization for this purpose.  Thus, we propose the following system over the domain $\Omega$.
\begin{subequations}\label{NSE_alpha}
 \begin{empheq}[left=\empheqlbrace]{align}
   \label{NSE_asml}
     (I-\alpha^2\Delta)\partial_t\bv
     +
     (\bv\cdot\nabla)\bv
     +
     \nabla q
     &=
     \nu\Delta\bv 
     + 
     \mathbf{f} 
     + 
     \mu (I_h(\bu) - I_h(\bv)),
     \\
   \label{NSE_alpha_div}
     \nabla\cdot \bv
     &= 
     0,
     \\
   \label{NSE_alpha0}
     \bv(\bx,0) &= \bv_0(\bx) = \bv_0.
 \end{empheq}
\end{subequations}
We take $I_h$ to be a linear operator satisfying
\begin{equation}
     \|\phi - I_{h}(\phi)\|_{L^2} \leq c_{1}h\|\nabla\phi\|_{L^2}
  \label{Asml}
\end{equation}
for some constants $c_1, h>0$, and all $\phi\in H^1$. Many common interpolants satisfy this bound.  For example, piecewise linear interpolation on a mesh with minimum spacing $h$ between grid points (see, e.g., \cite{Ern_Guermond_book1}, Corollary 1.109), and low-mode Fourier truncation preserving a number of modes proportional to $1/h$ (see, e.g., \cite{Canuto_Hussaini_Quarteroni_Zang_2006,Shen_Tang_Wang_2011}), both satisfy \eqref{Asml}.

We now state our global well-posedness result regarding system \eqref{NSE_alpha} and provide in Section~\ref{secProof1} a formal proof (using {\it a priori} estimates) of global existence of solution in $H$ and $V$. 
\begin{theorem}\label{T1}
Let $\bv_0 \in V$, and $\mathbf{f}\in H^1$ be time-independent. Then, for $\alpha, h, \mu>0$, system \eqref{NSE_alpha} possesses a unique solution $\bv(\bx, t)\in L^{\infty}(0, T; V)$ for any $T>0$. 
\end{theorem}

We recall the dimensionless Grashof number $G$, defined as 
$$G := \frac{\|\mathbf{f}\|_{L^2}}{\nu^2\lambda_1},$$ 
where $\lambda_1$ is the first eigenvalue of the Stokes operator $\bA$ defined in Section~\ref{secPre}. 

The next two theorems are our main results concerning the convergence of the data assimilation algorithm. 
\begin{theorem}\label{T2}
Suppose $\bu_0, \bv_0\in H^2\cap V$, and $\mathbf{f}\in V$. Let $\bu$ and $\bv$ be the corresponding solutions of \eqref{NSE} and \eqref{NSE_alpha}, respectively.  
Let $h>0$ and $\mu>0$ be such that 
$$h<\frac{\lambda_{1}^{-1/2}}{2c_{1}\sqrt{2C}}G^{-1} \text{\quad and \quad} \frac{\mu}{2} - C\nu\lambda_{1}G^{2} =: M_1 > 0,$$
for some constant $C>0$.
Then, for some $t_0>0$, and a.e. $t>t_0$ and any $\alpha>0$ such that 
\begin{equation}\label{alpha_bound_M1}
\alpha^2 < \frac{\nu}{M_1},
\end{equation}
it holds that
\begin{align}
 &\quad\|\bu(t) - \bv(t)\|_{L^2}^2 + \alpha\|\nabla\bu(t) - \nabla\bv(t)\|_{L^2}^2
\\&\leq\notag
C_{\bu, \mathbf{f}}\alpha^4(C_{M_{1}} - e^{-C_{M_{1}}t}) 
\\&\qquad\notag
+ C_{M_{1}}\left(\|\bu_0 - \bv_0\|_{L^2}^2 + \alpha\|\nabla\bu_0 - \nabla\bv_0\|_{L^2}^2\right)e^{-C_{M_{1}}t}.
\end{align}
In particular,
\begin{align*}
 \lim_{t\maps\infty}\|\bu(t) - \bv(t)\|_{L^2} \leq C\alpha^2,
\end{align*}
and 
\begin{align*}
 \lim_{t\maps\infty}\|\nabla\bu(t) - \nabla\bv(t)\|_{L^2} \leq C\alpha.
\end{align*}
\end{theorem}

Theorem \ref{T2} says that, for large times, the regularized solution approximates the true solution to the 2D NS equations, up to some error, and this error goes to zero as a power of $\alpha$ as $\alpha\rightarrow0$.  The next theorem provides faster rate of convergence for the $H^1$ norm, and also convergence in the $H^2$ norm, provided the data is smoother, and the nudging parameter and interpolation parameters, $\mu$ and $h$, satisfies stricter inequalities.

\begin{theorem}\label{T3}
Suppose $\bu_0, \bv_0\in H^3\cap V$, and $\mathbf{f}\in H^2\cap V$. 
Let $h, \mu>0$ be such that 
$$h<\frac{\lambda_{1}^{-1/2}}{2c_{1}\sqrt{2C}}G^{-1} \text{\quad and \quad} \frac{\mu}{2} - C\nu\lambda_{1}G^{8/3} =:M_2 > 0, $$
for some constant $C>0$ depending on $\nu$, $G$ and $\|\mathbf{f}\|_{H^2}$.
Then, for some $t_0>0$, for a.e. $t>0$ and any $\alpha>0$ such that 
$$\alpha^2 < \frac{\nu}{M_2},$$
the following inequality holds.
\begin{align}
 &\quad\|\nabla\bu(t) - \nabla\bv(t)\|_{L^2}^2 + \alpha\|\Delta\bu(t) - \Delta\bv(t)\|_{L^2}^2
\\&\leq\notag
\tilde{C}_{\bu, \mathbf{f}}\alpha^4(\tilde{C}_{M_{2}} - e^{-\tilde{C}_{M_{2}}t}) 
\\&\qquad\qquad\notag
+ \tilde{C}_{M_{2}}\left(\|\nabla\bu_0 - \nabla\bv_0\|_{L^2}^2 + \alpha\|\Delta\bu_0 - \Delta\bv_0\|_{L^2}^2\right)e^{-\tilde{C}_{M_{2}}t}.
\end{align}
In particular,
\begin{align*}
 \lim_{t\maps\infty}\|\nabla\bu(t) - \nabla\bv(t)\|_{L^2} \leq \tilde{C}\alpha^2,
\end{align*}
and 
\begin{align*}
 \lim_{t\maps\infty}\|\Delta\bu(t) - \Delta\bv(t)\|_{L^2} \leq \tilde{C}\alpha.
\end{align*}
\end{theorem}


\noindent
\section{Preliminaries}\label{secPre}
\noindent

All through this paper we denote the usual Lebesgue and Sobolev spaces by $L^{p}$ for $1\leq p\leq \infty$ and $H^{s}\equiv W^{s, 2}$ for $s > 0$, respectively. 
Let $\mathcal{V}$ be the set of all divergence-free, mean-free, smooth vector fields on $\mathbb{T}^2$. 
We follow the standard convention of denoting by $H$ and $V$ the closures of $\mathcal{V}$ in $L^2$ and $H^1$, respectively, 
with inner products 
$$(\bu, \bv) = \sum_{i=1}^2\int_{\mathbb{T}^2}\bu_{i}\bv_{i}\,d\bx \text{ \,\,and\,\,  } (\nabla \bu, \nabla \bv) = \sum_{i, j=1}^2\int_{\mathbb{T}^2}\partial_{j}\bu_{i}\partial_{j}\bv_{i}\,d\bx,$$
respectively, associated with the norms $\|\bu\|_{H}=(\bu, \bu)^{1/2}$ and $\| \bu \|_{V}=(\nabla \bu, \nabla \bu)^{1/2}$. For the sake of convenience, we use $\|\bu\|_{L^2}$ and $\|\bu\|_{H^1}$ to denote the above norms in $H$ and $V$, respectively. 

The following materials are standard in the study of fluid dynamics, 
in particular for the Navier-Stokes equations and related PDEs,  
and we refer to reader to \cite{Constantin_Foias_1988, Temam_2001_Th_Num} for more details. 
We define the Stokes operator $\bA\triangleq -\bP_{\sigma}\Delta$ 
with domain $\mathcal{D}(\bA)\triangleq H^2\cap V$, 
where $\bP_{\sigma}$ is the orthogonal Leray-Helmholtz projector from $L^2$ to $H$. 
Notice that under periodic boundary conditions, 
we have $\bA = -\Delta \bP_{\sigma}$.
Moreover, the Stokes operator can be extended 
as a linear operator from $V$ to $V'$ as 
$$\left<\bA\bu, \bv\right> = (\nabla \bu, \nabla \bv) \text{  for all  } \bv\in V.$$
It is well-known that $\bA^{-1} : H \hookrightarrow \mathcal{D}(\bA)$ 
is a positive-definite, self-adjoint, and compact operator from $H$ into itself,
thus, $\bA^{-1}$ possesses an orthonormal basis of positive eigenfunctions $\{ w_{k}\}_{k=1}^{\infty}$ in $H$, corresponding to a sequence of non-increasing sequence of eigenvalues. 
Therefore, $\bA$ has non-decreasing eigenvalues $\lambda_{k}$, 
i.e., $0 \leq \lambda_1 \leq \lambda_2, \ldots$. 
Then, we have the following Poincar\'e inqualities:
$$\lambda_1\|\bu\|_{L^2}^2\leq\|\nabla\bu\|_{L^2}^2 \text{\quad for\quad} \bu\in V,$$
$$\lambda_1\|\nabla\bu\|_{L^2}^2\leq\|\bA\bu\|_{L^2}^2 \text{\quad for\quad} \bu\in D(\bA).$$
Thus, $\|\nabla\bu\|_{L^2}$ is equivalent to $\|\bu\|_{H^1}$. 

Also, we state the following Ladyzhenskaya inequality 
$$\|\bu\|_{L^4}^2 \leq c\|\bu\|_{L^2} \|\nabla\bu\|_{L^2}$$
for all $\bu\in V$, which is a variation of the following interpolation result 
that is frequently used in this paper (see, e.g., \cite{Nirenberg_1959_AnnPisa} for a detailed proof).
Assume $1 \leq q, r \leq \infty$, and $0<\gamma<1$.  
For $v\in L^q(\mathbb{T}^{n})$, such that  $\partial^\alpha v\in L^{r} (\mathbb{T}^{n})$, for $|\alpha|=m$, then 
\begin{align}\label{PT1}
\|\partial_{s}v\|_{L^{p}} \leq C\|\partial^{\alpha}v\|_{L^{r}}^{\gamma}\| v\|_{L^{q}}^{1-\gamma},
\quad\text{where}\quad
\frac{1}{p} - \frac{s}{n} = \left(\frac{1}{r} - \frac{m}{n}\right) \gamma+ \frac{1}{q}(1-\gamma).
\end{align}

Next, for any $\bu, \bv \in \mathcal{V}$, 
we use the standard notation for the bilinear term
\begin{align*}
     \bB(\bu, \bv) := \bP_{\sigma}((\bu\cdot\nabla)\bv),
\end{align*}
which can be extended to a continuous map 
$\bB : V \times V \to V'$, where $V'$ denotes the dual space of $V_{\sigma}$, such that
\begin{align*}
     \left<\bB(\bu, \bv), \bw\right> = \int_{\mathbb{T}^2}(\bu\cdot\nabla \bv)\cdot \bw\,d\bx.
\end{align*}
for smooth functions $\bu, \bv, \bw\in \mathcal{V}$. 
 
We will use the following important properties of the map $B$, also known as Ladyzhenskaya-Sobolev inequalities which can be derived from the aforementioned Nirenberg inequality. Detailed proofs can be found in, e.g., \cite{Constantin_Foias_1988, Foias_Manley_Rosa_Temam_2001}.
\begin{lemma}
\label{L1}
For the operator $\bB$, the following identities hold.
\begin{subequations}
\begin{align}
    \label{symm1}
    \ip{\bB(\bu,\bv)}{\bw}_{V'} &= -\ip{\bB(\bu,\bw)}{\bv}_{V'}, 
    \quad\forall\;\bu\in V, \bv\in V, \bw\in V,\\
    \label{symm2}
    \ip{\bB(\bu,\bv)}{\bv}_{V'} &= 0,
    \quad\forall\;\bu\in V, \bv\in V, \bw\in V.
    \end{align}
\end{subequations}
Moreover, 
    \begin{subequations}
\begin{align}
    \label{B:424}
    |\ip{\bB(\bu,\bv)}{\bw}_{V'}|
    &\leq C\| \bu\|_{L^2}^{1/2} \|\nabla \bu\|_{L^2}^{1/2} \|\nabla \bv\|_{L^2} \|\bw\|_{L^2}^{1/2} \|\nabla \bw\|_{L^2}^{1/2},
\\
    \label{B:442}
    |\ip{\bB(\bu,\bv)}{\bw}_{V'}|
    &\leq C\| \bu\|_{L^2}^{1/2} \|\nabla \bu\|_{L^2}^{1/2} \|\nabla \bv\|_{L^2}^{1/2} \|\bA \bv\|_{L^2}^{1/2} \|\bw\|_{L^2},
    \\
    \label{B:inf22}
    |\ip{\bB(\bu,\bv)}{\bw}_{V'}|
    &\leq C\| \bu\|_{L^2}^{1/2} \|\bA \bu\|_{L^2}^{1/2} \|\nabla \bv\|_{L^2} \|\bw\|_{L^2},
\end{align}
for all $\bu$, $\bv$, $\bw$ in the largest spaces $H$, $V$, or $D(A)$, for which the right-hand sides of the inequalities are finite.
\end{subequations}
\end{lemma}

Due to the periodic boundary condition of our domain $\Omega = \mathbb{T}^2$, 
we have the additional orthogonality property for all $\bw\in D(\bA)$: 
\begin{align}\label{enstrophy_miracle}
\ip{\bB(\bw, \bw)}{\bA\bw} = 0,
\end{align}
and the Jacobi identity   
\begin{align}\label{jacobi}
 \ip{\bB(\bu, \bw)}{\bA\bw} + \ip{\bB(\bw, \bu)}{\bA\bw} +\ip{\bB(\bw, \bw)}{\bA\bu} = 0.
\end{align}
See, e.g., \cite{Constantin_Foias_1988,Temam_2001_Th_Num,Foias_Manley_Rosa_Temam_2001} for details.

Regarding the pressure term, we recall a corollary of a deep result of deRham; namely, that for any distribution $\mathbf{F}$, the equality $\mathbf{F}=\nabla p$ holds for some distribution $p$ if and only if $\left<\mathbf{F}, \bw\right> = 0$ for all $\bw \in \mathcal {V}$. 
See \cite{Wang_1993} for an elementary proof of the corollary.

In order to prove Theorem~\ref{T1}, we summarize a series of well-known results about 2D NS system in the next theorem and refer the readers to \cite{Constantin_Foias_1988, Temam_2001_Th_Num, Dascaliuc_Foias_Jolly_2005, Dascaliuc_Foias_Jolly_2007, Dascaliuc_2007, Dascaliuc_2008, Dascaliuc_2010} for more details. 
\begin{theorem}\label{T4}
With initial data $\bu_0\in V$ and external forcing $\mathbf{f}\in H$, the 2D NS system \eqref{NSE} possess a unique global solution $\bu\in L^{\infty}(0, \infty; V)\cap L^2(0, \infty; D(\bA))$. 
In particular, for some $t_0$ and $t>t_0$, we have $$\|\bu(t)\|_{L^2}^2 \leq 2\nu^{2}G^2\quad\text{ and }\quad \|\nabla\bu(t)\|_{L^2}^2 \leq 2\nu^{2}\lambda_{1}G^2. $$
Moreover, if $\mathbf{f} \in L^{2}(0, T; H)$, 
then the $H^2$ norm of $\bu$ is uniformly bounded after certain time $t_0>0$, 
i.e., we have 
$$\|\bA\bu\|_{L^2}^2 \leq C\nu^2\lambda_{1}^{2}G^2(G^2+C_{\mathbf{f}}) \sim C\nu^2\lambda_{1}^{2}G^4.$$ 
If we further assume that $\bu_0\in H^{2}\cap V$ and $\mathbf{f} \in H^2$, then the $H^3$ norm of $\bu$ is uniformly bounded; namely,
$$\|\bA^{3/2}\bu\|_{L^2}^2 \leq C\nu^2\lambda_{1}^{3}G^2(G^4+C_{\mathbf{f}}G^2+C_{\mathbf{f}}G+C_{\mathbf{f}}) \sim C\nu^2\lambda_{1}^{3}G^6,$$
while the square of the $H^{4}$ norm of $\bu$ is uniformly integrable in time, 
where $C_{\mathbf{f}}$ depends on $\|\bA\mathbf{f}\|_{L^2}$, $\|\mathbf{f}\|_{L^2}$, and $\lambda_1$. 
\end{theorem}
Also, we use the following Brezis-Gallouet inequality in our paper. 
\begin{lemma}
\label{LBG}
Let $u\in H^2(\Omega)$ where $\Omega \subset \mathbb{R}^2$. Then, there exists a constant $C$ such that $$\|u\|_{L^{\infty}} \leq C\|u\|_{H^1}\left(1+\log{\frac{\|\Delta u\|_{L^2}}{\lambda_1\|u\|_{H^1}}}\right)^{1/2}$$
\end{lemma}
For the sake of completeness, we state a special case of the generalized Gr\"onwall inequality proved in \cite{Jones_Titi_1992} (see also \cite{Farhat_Lunasin_Titi_2017_Horizontal}). 
\begin{lemma}
\label{L2}
Suppose that $Y(t)$ is a locally integrable and absolutely continuous function that satisfies the following: 
$$\frac{d Y}{d t} + \xi(t) Y \leq \beta(t), \quad\text{ a.e. on } (0, \infty), $$
such that 
$$\liminf_{t \to \infty} \int_{t}^{t+\tau} \xi(s)\,ds \geq \gamma, \quad\quad\quad \limsup_{t \to \infty} \int_{t}^{t+\tau} \xi^{-}(s)\,ds < \infty, $$
and 
$$\lim_{t \to \infty} \int_{t}^{t+\tau} \beta^{+}(s)\,ds = 0, $$
for fixed $\tau > 0$, and $\gamma > 0$, 
where 
$\xi^{-} = \max\{-\xi, 0\}$ 
and 
$\beta^{+} = \max\{\beta, 0\}$. 
Then, $Y(0) \to 0$ at an exponential rate as $t \to \infty$.
\end{lemma}

\begin{remark}
\label{R1}
In fact we will use the above lemma with a slight variation, i.e., 
$$\lim_{t \to \infty} \int_{t}^{t+\tau} \beta^{+}(s)\,ds \leq C\alpha^4, $$
which was proved in \cite{CarlsonHudsonLarios_2018} by a straight-forward minor modification of the proof in \cite{Jones_Titi_1992}. 
\end{remark}

\noindent
\section{Global well-posedness of system \eqref{NSE_alpha}}\label{secProof1}
\noindent
In this section, we provide the proof of Theorem~\ref{T1}. Here, we work formally, but we note that  following formal derivation including the {\it a priori} estimates can be made rigorous by using a Galerkin approximation then passing to the limit in an appropriate sense under certain compactness conditions by using Aubin-Lions Lemma (c.f. \cite{Constantin_Foias_1988,Temam_2001_Th_Num}). Since such approach is standard for both 2D and 3D NS system as well as other turbulence models such NS-Voigt system and NS-$\alpha$ model, we omit the details here but refer the readers to \cite{Constantin_Foias_1988, Larios_Pei_2017_MHDB_PS, Temam_2001_Th_Num, Farhat_Lunasin_Titi_2017_Horizontal,Larios_Titi_2009,Cao_Lunasin_Titi_2006} and the reference therein for more detailed discussions. 

{\smallskip\noindent {\em Proof of Theorem~\ref{T1}.}}

Taking a formal inner-product of \eqref{NSE_asml} with $\bv$, integrating by parts, and employing \eqref{symm2}, we obtain
\begin{align}\label{Energy}
     &\quad
     \frac{1}{2}\frac{d}{d t}\left(\|\bv\|_{L^2}^2 + \alpha^2\|\nabla\bv\|_{L^2}^2\right)
     +
     \nu\|\nabla\bv\|_{L^2}^2
     \\&\notag=
     \int_{\Omega}\mathbf{f}\cdot\bv\,d\bx
     +
     \mu\int_{\Omega}I_{h}(\bu)\cdot\bv\,d\bx
     -
     \mu\int_{\Omega}I_{h}(\bv)\cdot\bv\,d\bx
\end{align} 
The first term on the right side of \eqref{Energy} is bounded by 
\begin{align*}
     \|\mathbf{f}\|_{L^2} \|\bv\|_{L^2}
     \leq
     \frac{2}{\nu}\|\mathbf{f}\|_{L^2}^2
     +
     \frac{\nu}{8}\|\nabla\bv\|_{L^2}^2.
\end{align*}

Using \eqref{Asml}, we estimate the second term on the right side of \eqref{Energy} by 
\begin{align*}
     \mu\int_{\Omega}I_{h}(\bu)\cdot\bv\,d\bx
     &
     \leq
     \mu\int_{\Omega}(|\bu - I_{h}(\bu)| + |\bu|)|\bv|\,d\bx
     \\&\leq
     \mu\|\bu - I_{h}(\bu)\|_{L^2}\|\bv\|_{L^2}
     +
     \mu\|\bu\|_{L^2}\|\bv\|_{L^2}
     \\&
     \leq
     \frac{\mu c_{1}h}{\sqrt{\lambda_1}}\|\bu\|_{H^1}\|\nabla\bv\|_{L^2}
     +
     \frac{\mu}{\sqrt{\lambda_1}}\|\bu\|_{L^2}\|\nabla\bv\|_{L^2}
     \\&
     \leq
     \frac{2\mu^2 c_{1}^2h^2}{\lambda_1 \nu}\|\bu\|_{H^1}^2
     +
     \frac{2\mu^2}{\lambda_1 \nu}\|\bu\|_{L^2}^2
     +
     \frac{\nu}{4}\|\nabla\bv\|_{L^2}^2,
\end{align*}
where we used the Cauchy-Schwarz, Poincar\'e, and Young inequalities. 

Regarding the last term on the right side of \eqref{Energy}, we use \eqref{Asml} to obtain
\begin{align*}
     -
     \mu\int_{\Omega}I_{h}(\bv)\cdot\bv\,d\bx
     &
     =
     \mu\int_{\Omega}(\bv - I_{h}(\bv))\bv\,d\bx
     -
     \mu\|\bv\|_{L^2}^2
     \\&
     \leq
     \mu\|\bv - I_{h}(\bv)\|_{L^2}\|\bv\|_{L^2}
     -
     \mu\|\bv\|_{L^2}^2
     \\&
     \leq c_{1}h\mu\|\bv\|_{L^2}\|\nabla\bv\|_{L^2}
     -
     \mu\|\bv\|_{L^2}^2
     \\&
     \leq
     \left(\frac{2}{\nu} c_{1}^2h^2\mu^2 - \frac{\mu}{\nu}\right)\|\bv\|_{L^2}^2
     +
     \frac{\nu}{8}\|\nabla\bv\|_{L^2}^2. 
\end{align*}
Then, by combining all the above estimates, we obtain 
\begin{align*}
     &
     \frac{d}{d t}\left(\|\bv\|_{L^2}^2 + \alpha^2\|\nabla\bv\|_{L^2}^2\right)
     +
     \nu\|\nabla\bv\|_{L^2}^2
     \leq
     C
     +
     \left(\frac{2}{\nu} c_{1}^2h^2\mu^2 - \frac{\mu}{\nu}\right)\|\bv\|_{L^2}^2, 
\end{align*}
where the constant $C$ depends on $\lambda_1$, $\nu$, $\|\bu\|_{H^1}$, $\|\mathbf{f}\|_{L^2}$, $c_1$, $h$, and $\mu$. 
Therefore, Gr\"onwall inequality implies
\[\bv\in L^{\infty}(0, T; V)\]
for all $T>0$. 

In order to obtain the $H^2$ bounds on $\bv$, 
we multiply \eqref{NSE_asml} by $\bA\bv$ and integrate by parts to obtain
\begin{align}
      &\quad\frac{1}{2}\frac{d}{d t}\left(\|\nabla\bv\|_{L^2}^2 + \alpha^2\|\bA\bv\|_{L^2}^2\right)
     +
     \nu\|\bA\bv\|_{L^2}^2
     \\&=
     -
     \int_{\Omega}\mathbf{f}\cdot\bA\bv\,d\bx
     -
     \mu\int_{\Omega}I_{h}(\bu)\cdot\bA\bv\,d\bx
     +
     \mu\int_{\Omega}I_{h}(\bv)\cdot\bA\bv\,d\bx, 
     \label{Energy1}
\end{align}
where we used \eqref{symm2}. 
Then, we estimate the three terms on the right side of \eqref{Energy1}. 
Using Cauchy-Schwarz inequality, the first term is bounded by 
\begin{align*}
     \int_{\Omega}|\mathbf{f}|\,|\bA\bv|\,d\bx
     &
     \leq
     \|\mathbf{f}\|_{L^2} \|\bA\bv\|_{L^2}
     \leq
     \frac{2}{\nu}\|\mathbf{f}\|_{L^2}^2
     +
     \frac{\nu}{8}\|\bA\bv\|_{L^2}^2. 
\end{align*} 
By \eqref{Asml}, the second term is bounded by 
\begin{align*}
     \mu\int_{\Omega}(|\bu - I_{h}(\bu)| + |\bu|)|\bA\bv|\,d\bx
     &
     \leq
     \mu\|\bu - I_{h}(\bu)\|_{L^2}\|\bA\bv\|_{L^2}
     +
     \mu\|\bu\|_{L^2}\|\bA\bv\|_{L^2}
     \\&
     \leq
     \mu c_{1}h\|\bu\|_{H^1}\|\bA\bv\|_{L^2}
     +
     \mu\|\bu\|_{L^2}\|\bA\bv\|_{L^2}
     \\&
     \leq
     \frac{2\mu^2 c_{1}^2h^2}{\nu}\|\bu\|_{H^1}^2
     +
     \frac{2\mu^2}{\nu}\|\bu\|_{L^2}^2
     +
     \frac{\nu}{4}\|\bA\bv\|_{L^2}^2,    
\end{align*} 
while the last term is estimated similarly as 
\begin{align*}
     \mu\int_{\Omega}I_{h}(\bv)\cdot\bA\bv\,d\bx
     &
     =
     -
     \mu\int_{\Omega}(\bv - I_{h}(\bv))\bA\bv\,d\bx
     +
     \mu\int_{\Omega}\bv\cdot\bA\bv\,d\bx
     \\&
     \leq
     \mu\|\bv - I_{h}(\bv)\|_{L^2}\|\bA\bv\|_{L^2}
     -
     \mu\|\nabla\bv\|_{L^2}^2
     \\&
     \leq
     c_{1}h\mu\|\nabla\bv\|_{L^2}\|\bA\bv\|_{L^2}
     -
     \mu\|\nabla\bv\|_{L^2}^2
     \\&
     \leq
     \left( \frac{2c_{1}^2 h^2\mu^2}{\nu} - \mu \right)\|\nabla\bv\|_{L^2}^2
     +
     \frac{\nu}{8}\|\bA\bv\|_{L^2}^2. 
\end{align*} 
Combining the above estimates, we obtain
\begin{align*}
     \frac{d}{d t} \left( \|\nabla\bv\|_{L^2}^2 + \alpha^2\|\bA\bv\|_{L^2}^2\right)
     +
     \nu\|\bA\bv\|_{L^2}^2
     \leq
     C
     +
     \left( \frac{2c_{1}^2 h^2\mu^2}{\nu} - \mu \right)\|\nabla\bv\|_{L^2}^2,
\end{align*}
where the constant $C$ depends on $\lambda_1$, $\nu$, $\|\bu\|_{H^1}$, $\|\mathbf{f}\|_{L^2}$, $c_1$, $h$, and $\mu$. 
Hence, Gr\"onwall inequality implies $\bv \in L^{\infty}(0, T; H^2\cap V)$. 

Next, we show that the time derivative $d\bv/dt$ is uniformly bounded in $V$. 
In fact, for any test function $\phi\in V$, we have 
\begin{align}
     \|(I - \alpha^2\Delta)\frac{d\bv}{dt}\|_{V'} 
     &\leq
     \sup_{\|\phi\|_{V}=1}\ip{\bB(\bv, \bv)}{\phi}
     +
     \sup_{\|\phi\|_{V}=1}\nu\ip{\bA\bv}{\phi}
     +
     \sup_{\|\phi\|_{V}=1}\ip{\mathbf{f}}{\phi} 
     \label{Ineq}
     \\&\quad
     +
     \sup_{\|\phi\|_{V}=1}\mu\ip{I_{h}(\bu)}{\phi}
     +
     \sup_{\|\phi\|_{V}=1}\mu\ip{I_{h}(\bv)}{\phi}, \nonumber
\end{align}
where we integrated by parts and used $\nabla\cdot\phi = 0$. 
Since $\bv$ is bounded in $V$, 
the first term on the right side of \eqref{Ineq} is bounded by 
\begin{align*}
     &
     C\|\bv\|_{L^3}\|\nabla\bv\|_{L^2}\|\phi\|_{L^6}
     \leq
     C\|\bv\|_{L^2}^{1/2}\|\nabla\bv\|_{L^2}^{3/2}\|\nabla\phi\|_{L^2}
     \leq
     C\|\bv\|_{L^2}^{1/2}\|\bv\|_{V}.
\end{align*}
After integration by parts, the second term is bounded by 
\begin{align*}
     \nu\|\nabla\bv\|_{L^2}\|\nabla\phi\|_{L^2}
     \leq
     \nu\|\bv\|_{V},
\end{align*}
while Cauchy-Schwarz inequality implies that the third term is bounded by 
$$\|\mathbf{f}\|_{L^2}\|\phi\|_{L^2} \leq C\|\mathbf{f}\|_{L^2}. $$
Regarding the last two terms, we used \eqref{Asml} and proceed as 
\begin{align*}
     \sup_{\|\phi\|_{V}=1}\mu\ip{I_{h}(\bu)}{\phi}
     &\leq
     \sup_{\|\phi\|_{V}=1}\mu\ip{I_{h}(\bu) - \bv}{\phi}
     +
     \sup_{\|\phi\|_{V}=1}\mu\ip{\bu}{\phi}
     \\&
     \leq
     \mu\|I_{h}(\bu) - \bu\|_{L^2}\|\phi\|_{L^2}
     +
     \mu\|\bu\|_{L^2}\|\phi\|_{L^2}
     \\&
     \leq
     c_{1}h\mu\|\nabla\bu\|_{L^2}\|\phi\|_{L^2}
     +
     \mu\|\bu\|_{L^2}\|\phi\|_{L^2}
     \\&
     \leq
     C\mu\|\bu\|_{V},
\end{align*}
and similar argument gives 
\begin{align*}
     &
     \sup_{\|\phi\|_{V}=1}\mu\ip{I_{h}(\bv)}{\phi}
     \leq
     C\mu\|\bv\|_{V}. 
\end{align*}
Thus, we conclude that $(I-\alpha^2\Delta)d\bv/dt$ is uniformly bounded in $V'$. 
By inverting the Helmholtz operator $(I-\alpha^2\Delta)$, 
we obtain the uniform boundedness of $d\bv/dt$ in $V$. 
By using a difference test function $\tilde{\phi} \in D(\bA)$, 
similar estimates provide that $(I-\alpha^2\Delta)d\bv/dt$ is uniformly bounded in $D(\bA)'$, 
which in turn indicates that $d\bv/dt$ is uniformly bounded in $H$. 
Therefore, together with the $H^2$ bounds on $\bv$, 
one can apply standard compactness arguments and pass to the limit and get the existence of a solution in $V$. 

The uniqueness is proven in a similar way to standard arguments for Navier-Stokes system, see, e.g.,  \cite{Constantin_Foias_1988,Temam_2001_Th_Num}. Hence, for the sake of simplicity, we omit the details here. 
The proof of Theorem~\ref{T1} is thus complete. 

\noindent
\section{Proof of $L^2$ convergence rate}\label{secProof2}
\noindent
In this section, we provide the proof of Theorem~\ref{T2}. 

{\smallskip\noindent {\em Proof of Theorem~\ref{T2}.}}

By denoting the difference $\bw = \bu - \bv$, $\Pi = p -q$ and taking the difference of \eqref{NSE1} and \eqref{NSE_asml}, we get 
\begin{subequations}\label{Diff}
 \begin{empheq}[left=\empheqlbrace]{align}
     \partial_{t}\bw
     -
     \nu\Delta\bw
     +
     \alpha^2\Delta\bu_{t}
     -
     \alpha^2\Delta\bw_{t}
     &=
     -
     (\bw\cdot\nabla)\bu
     -
     (\bu\cdot\nabla)\bw \nonumber
     \\&\quad
     +
     (\bw\cdot\nabla)\bw
     -
     \nabla\Pi
     -\mu I_{h}(\bw), 
     \label{Diff_eq}
     \\
     \nabla\cdot \bw
     &= 
     0,
     \\
     \bw(\bx,0) &= \bu_0(\bx) - \bv_0(\bx).
 \end{empheq}
\end{subequations}
Multiply \eqref{Diff_eq} by $\bw$, integrate by parts over $\Omega$, and we obtain 
\begin{align}
     &
     \frac{1}{2}\frac{d}{d t}\left(\|\bw\|_{L^2}^2 + \alpha^2\|\nabla\bw\|_{L^2}^2\right)
     +
     \nu\|\nabla\bw\|_{L^2}^2 \nonumber
     \\&\quad
     =
     -
     \alpha^2\int_{\Omega}\Delta(\partial_{t}\bu)\cdot\bw\,d\bx
     -
     \int_{\Omega}(\bw\cdot\nabla)\bu\cdot\bw\,d\bx
     -
     \mu\int_{\Omega}I_{h}(\bw)\cdot\bw\,d\bx,
   \label{Energy2}
\end{align}
where we used \eqref{symm2}. 
For the first term on the right side of \eqref{Energy2}, 
we integrate by parts and use \eqref{NSE1}in order to get 
\begin{align}
     -\alpha^2\int_{\Omega}\Delta(\partial_{t}\bu)\cdot\bw\,d\bx
     &=
     -\alpha^2\int_{\Omega}\partial_{t}\bu\cdot\Delta\bw\,d\bx \nonumber
     \\&
     =
     -
     \alpha^2\int_{\Omega}\nu\Delta\bu\cdot\Delta\bw\,d\bx
     +
     \alpha^2\int_{\Omega}(\bu\cdot\nabla)\bu\cdot\Delta\bw\,d\bx \nonumber
     \\&\quad\quad
     +
     \alpha^2\int_{\Omega}\nabla p\cdot\Delta\bw\,d\bx
     -
     \alpha^2\int_{\Omega}\mathbf{f}\cdot\Delta\bw\,d\bx.
   \label{Est}
\end{align}
Then, we estimate the four integrals on the right side of \eqref{Est}. 

Regarding the first integral, we integrate by parts and proceed as 
\begin{align*}
     -
     \alpha^2\nu\int_{\Omega}\Delta\bu\cdot\Delta\bw\,d\bx
     &\leq
     \alpha^2\nu\int_{\Omega}|\nabla\Delta\bu|\,|\nabla\bw|\,d\bx
     \leq
     \alpha^2\nu\|\nabla\Delta\bu\|_{L^2}\|\nabla\bw\|_{L^2}
     \\&
     \leq
     4\alpha^4\nu\|\bu\|_{H^3}^2
     +
     \frac{\nu}{16}\|\nabla\bw\|_{L^2}^2.
\end{align*}
After integration by parts, we bound the second integral by 
\begin{align*}
     &
     \alpha^2\int_{\Omega}|\nabla\bu|^2\,|\nabla\bw|\,d\bx
     +
     \alpha^2\int_{\Omega}|\bu|\,|\nabla\nabla\bu|\,|\nabla\bw|\,d\bx
     \\&\quad\quad
     \leq
     C\alpha^2\|\nabla\bu\|_{L^4}^2\|\nabla\bw\|_{L^2}
     +
     C\alpha^2\|\bu\|_{L^{\infty}}\|\Delta\bu\|_{L^2}\|\nabla\bw\|_{L^2}
     \\&\quad\quad
     \leq
     C\alpha^2\|\bu\|_{H^1}\|\bu\|_{H^2}\|\nabla\bw\|_{L^2}
     \\&\quad\quad\quad
     +
     C\alpha^2\|\bu\|_{H^1}\|\Delta\bu\|_{L^2}\left(1+\log{\frac{\|\Delta\bu\|_{L^2}}{\lambda_1\|\bu\|_{H^1}}}\right)^{1/2}\|\nabla\bw\|_{L^2}   
     \\&\quad\quad
     \leq
     \frac{C\alpha^4}{\nu}\|\bu\|_{H^1}^{2}\|\bu\|_{H^2}^{2}
     +
     \frac{C\alpha^4}{\nu}\|\bu\|_{H^1}^2\|\bu\|_{H^2}^2\left(1+\log{\frac{\|\Delta\bu\|_{L^2}}{\lambda_1\|\bu\|_{H^1}}}\right)
     \\&\quad\quad\quad
     +
     \frac{\nu}{16}\|\nabla\bw\|_{L^2}^2,
\end{align*}
where we used Lemma~\ref{L1} and Lemma~\ref{LBG}. 

The third integral is zero due to $\nabla\cdot\bw = 0$.

As for the last integral, 
we integrate by parts and apply Cauchy-Schwarz inequality and obtain 
\begin{align*}
     -
     \alpha^2\int_{\Omega}\mathbf{f}\cdot\Delta\bw\,d\bx
     &\leq
     \alpha^2\int_{\Omega}|\nabla\mathbf{f}|\,|\nabla\bw|\,d\bx
     \leq
     C\alpha^2\|\nabla\mathbf{f}\|_{L^2} \|\nabla\bw\|_{L^2}
     \\&
     \leq
     \frac{C\alpha^4}{\nu}\|\mathbf{f}\|_{H^1}^2
     +
     \frac{\nu}{16}\|\nabla\bw\|_{L^2}^2.
\end{align*}

Next, using H\"older inequality, 
we bound the second term on the right side of \eqref{Energy2} by 
\begin{align*}
     C\|\bw\|_{L^4} \|\nabla\bu\|_{L^2} \|\bw\|_{L^4}
     &\leq
     C\|\bu\|_{H^1}\|\bw\|_{L^2}\|\nabla\bw\|_{L^2}
     \\&
     \leq
     \frac{C}{\nu}\|\bu\|_{H^1}^{2}\|\bw\|_{L^2}^{2}
     +
     \frac{\nu}{16}\|\nabla\bw\|_{L^2}^2,
\end{align*}
where we used Lemma~\ref{L1} and Young's inequality. 

Finally, the last term of \eqref{Energy2} is estimated as 
\begin{align*}
     -
     \mu\int_{\Omega}I_{h}(\bw)\cdot\bw\,d\bx
     &=
     \mu\int_{\Omega} \bigl(\bw - I_{h}(\bw)\bigr)\cdot\bw\,d\bx
     -
     \mu\|\bw\|_{L^2}^2
     \\&
     \leq
     \mu\|\bw - I_{h}(\bw)\|_{L^2}\|\bw\|_{L^2}
     -
     \mu\|\bw\|_{L^2}^2
     \\&
     \leq
     c_{1}h\mu\|\nabla\bw\|_{L^2}\|\bw\|_{L^2}
     -
     \mu\|\bw\|_{L^2}^2
     \\&
     \leq
     \frac{2c_{1}^{2}h^2\mu^2}{\nu}\|\bw\|_{L^2}
     +
     \frac{\nu}{8}\|\nabla\bw\|_{L^2}^2
     -
     \mu\|\bw\|_{L^2}^2. 
\end{align*}
Combining all the above estimates and denoting
$$X(t) := \|\bw(t)\|_{L^2}^2+\alpha^2\|\nabla\bw\|_{L^2}^2 \text{\quad and \quad} Y(t) := \|\nabla\bw(t)\|_{L^2}^2,$$
we obtain 
\begin{align*}
     \frac{d}{d t}X
     +
     \nu Y
     &\leq
     C_{\bu,\mathbf{f}}\alpha^4
     +
     \left( \frac{C}{\nu}\|\bu\|_{H^1}^2 + \frac{2c_{1}^2h^2\mu^2}{\nu} - \mu\right)\|\bw\|_{L^2}^2,
\end{align*}
which is equivalent to 
\begin{align}
     \frac{d}{d t}X
     +
     \nu Y
     +
     \left( \mu - \frac{C}{\nu}\|\bu\|_{H^1}^2 - \frac{2c_{1}^{2}h^2\mu^2}{\nu} \right) \|\bw\|_{L^2}^2
     &\leq
     C_{\bu,\mathbf{f}}\alpha^4
  \label{Ineq1} 
\end{align}
where the constant $C_{\bu, \mathbf{f}}$ depends on $\nu$, $\|\bu\|_{H^3}$, $\|\bu\|_{H^2}$ and $\|\mathbf{f}\|_{H^1}$. 
Then, by our choice of $h$, we have 
$$\frac{2c_{1}^2h^2\mu^2}{\nu} < \frac{\mu}{2}.$$
By the fact that $\|\bw\|_{L^2}^2 = X - \alpha^{2}Y$, 
we obtain 
\begin{align*}
     \frac{d}{d t}X
     +
     \left( \frac{\mu}{2} - \frac{C}{\nu}\|\bu\|_{H^1}^2 \right) X
     +
     \left( \nu - \alpha^2 \left(\frac{\mu}{2} - \frac{C}{\nu}\|\bu\|_{H^1}^2\right) \right) Y
     &\leq
     C_{\bu,\mathbf{f}}\alpha^4 
\end{align*} 
Thus, in view of Theorem~\ref{T4}, for $t > t_0$, $$\nu - \alpha^2 \left(\frac{\mu}{2} - \frac{C}{\nu}\|\bu\|_{H^1}^2\right) > 0$$ due to our choice of $\alpha$ and $\mu$, from which we also infer that $\xi(t):=\frac{\mu}{2} - \frac{C}{\nu}\|\bu\|_{H^1}^2$ satisfies the condition in Lemma~\ref{L2}. 
Therefore, we get 
$$X(t) \leq C_{\bu, \mathbf{f}}\alpha^4(C_{M_{1}} - e^{-C_{M_{1}}t}) + C_{M_{1}}X(0)e^{-C_{M_{1}}t}$$
for all $t\geq t_0$.
Thus, proof of Theorem~\ref{T1} is complete. 

\noindent
\section{Proof of $H^1$ convergence rate}\label{secProof3}
\noindent
In this section, we provide the proof of Theorem~\ref{T3}. 

{\smallskip\noindent {\em Proof of Theorem~\ref{T3}.}}

Multiply equation \eqref{Diff_eq} by $\bA\bw$, 
integrate by parts over $\Omega$, and we get 
\begin{align}
     &
     \frac{1}{2}\frac{d}{d t}\left(\|\nabla\bw\|_{L^2}^2 + \alpha^2\|\Delta\bw\|_{L^2}^2\right)
     +
     \nu\|\Delta\bw\|_{L^2}^2 \nonumber
     \\&\quad\quad
     =
     \alpha^2\int_{\Omega}\Delta\bu_{t}\cdot\Delta\bw\,d\bx
     -
     \int_{\Omega}(\bw\cdot\nabla)\bw\cdot\Delta\bu\,d\bx
     +
     \mu\int_{\Omega}I_{h}(\bw)\cdot\Delta\bw\,d\bx,
  \label{Enstrophy2}
\end{align}
where we used \eqref{enstrophy_miracle} and \eqref{jacobi}.
In order to bound the first term on the right side of \eqref{Enstrophy2}, 
we integrate by parts and use the equation of $\bu$, and obtain 
\begin{align}
     \alpha^2\int_{\Omega}\Delta\bu_{t}\cdot\Delta\bw\,d\bx
     &=\notag
     \alpha^2\int_{\Omega}\bu_{t}\cdot\Delta^2\bw\,d\bx
     \\&=\label{Eq1}
     \alpha^2\int_{\Omega}\nu\Delta\bu\cdot\Delta^2\bw\,d\bx
     -
     \alpha^2\int_{\Omega}(\bu\cdot\nabla)\bu\cdot\Delta^2\bw\,d\bx 
     \\&\quad\quad\notag
     -
     \alpha^2\int_{\Omega}\nabla p\cdot\Delta^2\bw\,d\bx
     +
     \alpha^2\int_{\Omega}\mathbf{f}\cdot\Delta^2\bw\,d\bx.
\end{align}
Next, we estimate the four integrals on the right side of \eqref{Eq1}. 
After integration by parts, 
the first one is bounded by 
\begin{align*}
     \alpha^2\nu\|\Delta^2\bu\|_{L^2} \|\Delta\bw\|_{L^2}
     &
     \leq
     \alpha^4\nu\|\bu\|_{H^4}^2
     +
     \frac{\nu}{8} \|\Delta\bw\|_{L^2}^{2}.
\end{align*}
For the second integral, 
using integration by parts and Lemma~\ref{L1}, 
we have
\begin{align*}
     &
     -
     \alpha^2\int_{\Omega}(\bu\cdot\nabla)\bu\cdot\Delta^2\bw\,d\bx
     \\&\quad
     \leq
     3\alpha^2\int_{\Omega}|\nabla\bu|\,|\nabla^2\bu|\,|\Delta\bw|\,d\bx
     +
     \alpha^2\int_{\Omega}|\bu|\,|\nabla^3\bu|\,|\Delta\bw|\,d\bx
     \\&\quad
     \leq
     C\alpha^2\|\nabla\bu\|_{L^4} \|\nabla^2\bu\|_{L^4} \|\Delta\bw\|_{L^2}
     +
     C\alpha^2\|\bu\|_{L^{\infty}} \|\nabla^3\bu\|_{L^2} \|\Delta\bw\|_{L^2}
     \\&\quad
     \leq
     \frac{C\alpha^4}{\nu}\|\bu\|_{H^2}^2\|\bu\|_{H^3}^2
     +
     \frac{C\alpha^4}{\nu}\|\bu\|_{H^1}^2\|\bu\|_{H^3}^2\left(1+\log{\frac{\|\Delta\bu\|_{L^2}}{\lambda_1\|\bu\|_{H^1}}}\right)
     \\&\quad\quad
     +
     \frac{\nu}{4} \|\Delta\bw\|_{L^2}^{2}, 
\end{align*}
where we also used Lemma~\ref{LBG}. 
The third integral vanished due to $\nabla\cdot\bw = 0$. 

As for the last integral, 
we integrate by parts and use H\"older inequality in order to get 
\begin{align*}
     &
     \alpha^2\int_{\Omega}\mathbf{f}\cdot\Delta^2\bw\,d\bx
     =
     \alpha^2\int_{\Omega}\Delta\mathbf{f}\,\Delta\bw\,d\bx
     \\&\quad
     \leq
     \alpha^2\|\Delta\mathbf{f}\|_{L^2}\|\Delta\bw\|_{L^2}
     \leq
     2\alpha^4\|\mathbf{f}\|_{H^2}^2
     +
     \frac{1}{8} \|\Delta\bw\|_{L^2}^{2}.
\end{align*}

The second term on the right side of \eqref{Enstrophy2} is estimated as 
\begin{align*}
     &
     -
     \int_{\Omega}(\bw\cdot\nabla)\bw\cdot\Delta\bu\,d\bx
     \leq
     C\|\bw\|_{L^4} \|\nabla\bw\|_{L^4} \|\Delta\bu\|_{L^2}
     \\&\quad
     \leq
     C\lambda_1^{-1/4}\|\bu\|_{H^2}\|\nabla\bw\|_{L^2}^{3/2} \|\Delta\bw\|_{L^2}^{1/2}
     \\&
     \leq
     \frac{C}{(\nu\lambda_1)^{1/3}}\|\bu\|_{H^2}^{4/3} \|\nabla\bw\|_{L^2}^{2}
     +
     \frac{\nu}{8} \|\Delta\bw\|_{L^2}^{2},
\end{align*} 
where we used Poincar\'e's inequality. 

Regarding the last term in \eqref{Enstrophy2}, we have 
\begin{align*}
     \mu\int_{\Omega}I_{h}(\bw)\cdot\Delta\bw\,d\bx
     &=
     -
     \mu\int_{\Omega}\bigl(\bw - I_{h}(\bw)\bigr)\,\Delta\bw\,d\bx
     +
     \mu\int_{\Omega}\bw\cdot\Delta\bw\,d\bx
     \\&
     \leq
     \mu\|\bw - I_{h}(\bw)\|_{L^2} \|\Delta\bw\|_{L^2}
     -
     \mu\|\nabla\bw\|_{L^2}^2
     \\&
     \leq
     \frac{2c_{1}^2h^2\mu^2}{\nu}\|\nabla\bw\|_{L^2}^2
     +
     \frac{\nu}{8} \|\Delta\bw\|_{L^2}^{2}
     -
     \mu\|\nabla\bw\|_{L^2}^2.
\end{align*} 
By summing up all the above estimates and denoting 
$$\tilde{X}(t) = \|\nabla\bw(t)\|_{L^2}^2+\alpha^2\|\Delta\bw\|_{L^2}^2 \text{\quad and \quad} \tilde{Y}(t) = \|\Delta\bw(t)\|_{L^2}^2, $$
we obtain 
\begin{align}
     \frac{d}{d t}\tilde{X}(t)
     +
     \tilde{Y}(t)
     &\leq
     \tilde{C}_{\bu,\mathbf{f}}\alpha^4
     +
     \left(\frac{C}{(\nu\lambda_1)^{1/3}}\|\bu\|_{H^2}^{4/3} + \frac{2c_{1}^2h^2\mu^2}{\nu} - \mu\right)\|\nabla\bw\|_{L^2}^2,
  \label{Ineq2} 
\end{align}
where the time integral of $\tilde{C}_{\bu, \mathbf{f}}$ is bounded since $\bu\in L^{2}(0, T;{H^4})$ and $\|\mathbf{f}\|_{H^2}$. 
Hence, by our choice of $h$, $\mu$, and $\alpha$, we have  
$$\frac{2c_{1}^2h^2\mu^2}{\nu} < \frac{\mu}{2}$$ 
and for $t > t_0$,
$$\nu - \alpha^2\left(\frac{\mu}{2} - \frac{C}{(\nu\lambda_1)^{1/3}}\|\bu\|_{H^2}^{4/3}\right) \geq \nu - M_{2}\alpha^2 > 0.$$ 
Thus, similar arguments together with Lemma~\ref{L2} imply  
$$\tilde{X}(t) \leq \tilde{C}_{\bu, \mathbf{f}}\alpha^4(\tilde{C}_{M_{2}} - e^{-\tilde{C}_{M_{2}}t}) + \tilde{C}_{M_{2}}X(0)e^{-\tilde{C}_{M_{2}}t}$$
for all $t\geq t_0$.
The proof of Theorem~\ref{T3} is thus complete. 

\begin{remark}
There seems to be a compromise between the choice of $\mu$ and $\alpha$. According to Theorem \ref{T2}, $\mu\alpha^2$ should not be too large. However, choosing larger $\alpha$ provides a stronger regularizing effect to the solution of \eqref{NSE_alpha}, and choosing $\mu$ larger increases the convergence rate of the data assimilation. In a future computational work, we will examine the interplay between these parameters, and their effect on the accuracy and efficiency of approximating the true solution.
\end{remark}

\section*{Acknowledgement}
\noindent
  The research of A.L. was supported in part by NSF grant number DMS-1716801.
\begin{scriptsize}

\end{scriptsize}

\end{document}